\documentclass[12pt,leqno]{article}
\usepackage{amsmath, amsthm, amsfonts, amssymb, color}
\usepackage{mathrsfs}
\newcommand{\be}{\begin{eqnarray}}
\newcommand{\ee}{\end{eqnarray}}
\newcommand{\ce}{\begin{eqnarray*}}
\newcommand{\de}{\end{eqnarray*}}

\theoremstyle{definition}

\definecolor{wco}{rgb}{0.5,0.2,0.3}

\numberwithin{equation}{section} \theoremstyle{remark}

\newcommand{\ua}{\uparrow}

\title{{\bf Harnack Inequality and Strong Feller Property for Stochastic Fast-Diffusion Equations}
\footnote{Supported in part by WIMCS,
     NNSFC(10121101),
RFDP(20040027009), the 973-Project in China and DFG--Internationales Graduiertenkolleg ``Stochastics and Real World Models''.} }
\author{{\bf Wei Liu$^{a,c}$ and  Feng-Yu Wang$^{a,b}$\footnote{Corresponding author: wangfy@bnu.edu.cn}}\\
{\footnotesize $a.$ School of Mathematics, Beijing Normal
University,
Beijing 100875, China}\\
\footnotesize{$b.$ Department of Mathematics, Swansea University, Singleton Park, SA2 8PP, UK}
\\
\footnotesize{$c.$ Fakult\"at F\"ur Mathematik, Universit\"at
Bielefeld, D-33501 Bielefeld, Germany}}

\begin{document}
\maketitle
\def\R{\mathbb R}  \def\ff{\frac} \def\ss{\sqrt} \def\BB{\mathbf
B}
\def\N{\mathbf N} \def\kk{\kappa} \def\m{{\bf m}}
\def\dd{\delta} \def\DD{\Delta} \def\vv{\varepsilon} \def\rr{\rho}
\def\<{\langle} \def\>{\rangle} \def\GG{\Gamma} \def\gg{\gamma}
  \def\nn{\nabla} \def\pp{\partial} \def\tt{\tilde}
\def\d{\text{\rm{d}}} \def\bb{\beta} \def\aa{\alpha} \def\D{\mathcal D}
\def\E{\mathbb E} \def\si{\sigma} \def\ess{\text{\rm{ess}}}
\def\beg{\begin} \def\beq{\begin{equation}}  \def\F{\mathcal F}
\def\Ric{\text{\rm{Ric}}} \def\Hess{\text{\rm{Hess}}}\def\B{\mathbf B}
\def\e{\text{\rm{e}}} \def\ua{\underline a} \def\OO{\Omega} \def\b{\mathbf b}
\def\oo{\omega}     \def\tt{\tilde} \def\Ric{\text{\rm{Ric}}}
\def\cut{\text{\rm{cut}}} \def\P{\mathbb P} \def\ifn{I_n(f^{\bigotimes n})}
\def\fff{f(x_1)\dots f(x_n)} \def\ifm{I_m(g^{\bigotimes m})} \def\ee{\varepsilon}
\def\pm{\pi_{{\bf m}}}   \def\p{\mathbf{p}}   \def\ml{\mathbf{L}}
 \def\C{\mathcal C}      \def\aaa{\mathbf{r}}     \def\r{r}
\def\gap{\text{\rm{gap}}} \def\prr{\pi_{{\bf m},\varrho}}  \def\r{\mathbf r}
\def\Z{\mathbb Z} \def\vrr{\varrho} \def\ll{\lambda}
\def\L{\mathcal L}

\def\[{{\Big[}}
\def\]{{\Big]}}

\def\({{\Big(}}
\def\){{\Big)}}

\def\xy{{\|X_t-Y_t\|}}
\def\x{{\|X_t\|}}
\def\y{{\|Y_t\|}}
\def\mxy{{\mathbf{m}(|X_t-Y_t|^2(|X_t|\vee|Y_t|)^{r-1})}}
\def\mx{{\(\mathbf{m}[(|X_t|\vee|Y_t|)^{r+1}]\)^{\frac{1-r}{1+r}}}}

\def\bt{\begin{theorem}}
\def\et{\end{theorem}}
\def\bl{\begin{lemma}}
\def\el{\end{lemma}}
\def\br{\begin{remark}}
\def\er{\end{remark}}
\def\bx{\begin{Example}}
\def\ex{\end{Example}}
\def\bd{\begin{definition}}
\def\ed{\end{definition}}
\def\bp{\begin{proposition}}
\def\ep{\end{proposition}}
\def\bc{\begin{corollary}}
\def\ec{\end{corollary}}

\begin{abstract} As a continuation to \cite{W07}, where the Harnack
inequality and the strong Feller property are studied for a class of
stochastic generalized porous media equations, this paper presents
analogous results for stochastic fast-diffusion equations. Since the
fast-diffusion equation possesses weaker dissipativity than the
porous medium one does, some technical difficulties appear in the
study. As a compensation to the weaker dissipativity condition, a
Sobolev-Nash inequality is assumed for the underlying self-adjoint
operator in applications. Some concrete examples are constructed to
illustrate the main results.
\end{abstract}

\bigbreak

\section{Introduction}

Recently, the dimension-free Harnack inequality introduced in
\cite{W97} was established in \cite{W07} for a class of  stochastic
generalized porous media equations. As applications, the strong
Feller property, estimates of the transition density and some
contractivity properties were obtained for the associated Markov
semigroup. The approach used in \cite{W07} is based on a coupling
argument developed in  \cite{ATW}, where Harnack inequalities are
studied for diffusion semigroups on Riemannian manifolds with
unbounded below curvatures. The advantage of this approach is that
it avoids the assumption on curvature lower bounds used in previous
articles (see \cite{AK,AZ,BGL,RW1,RW2}), which is very hard to
verify in the framework of non-linear SPDEs. The main aim of this
paper is to apply this method to stochastic generalized
fast-diffusion equations studied in \cite{RRW}.

Let $(E,\mathcal M,{\bf m})$ be a separable probability space and
$(L,\D(L))$ a negative definite self-adjoint linear operator on
$L^2({\bf m})$ having
 discrete spectrum. Let
$$
  (0<) \ll_1\le \ll_2\le \cdots
$$ be all eigenvalues of $-L$ with unit eigenfunctions $\{e_i\}_{i\ge 1}$.

Next, let $H$ be the completion of $L^2(\m)$ under the inner product
$$\<x,y\>_H:= \sum_{i=1}^\infty \ff 1 {\ll_i}
\<x,e_i\>\<y,e_i\>,$$ where and in what follows, $\<\cdot ,\cdot \>$
is the inner product in $L^2(\m)$.  Let $L_{HS}$ denote the space of
all Hilbert-Schmidt operators from $L^2(\m)$ to $H$. Let $W_t$ be
the cylindrical Brownian motion on $L^2(\m)$ w.r.t. a complete
filtered probability space $(\OO, \F_t, \P)$; that is, $W_t= \{B_t^i
e_i\}_{i\ge 1}$ for a sequence of independent one-dimensional
$\F_t$-Brownian motions $\{B_t^i\}_{i\ge 1}$.

Let
$$\Psi: [0,\infty)\times \R\times \OO\to \R$$
be progressively measurable and continuous in the second variable, and let
$$ Q: [0,\infty)\times \OO\to L_{HS}$$
be progressively measurable such that for all $t\geq0$
\beq\label{1.0}
\|Q_t\|_{L_{HS}}^2\leq q_t,\ a.s.
\end{equation}
for $q_t$ some deterministic local integrable function on
$[0,+\infty)$. We consider the equation \beq\label{1.1}
 \d X_t = \big\{L\Psi(t,X_t)+ \gg_t X_t\big\}\d
t +Q_t \d W_t,
\end{equation}
where $\gg: [0,\infty)\to \R$  is locally bounded and measurable. In
particular, if $\gg=0, Q=0$ and $\Psi(t,s)=s^r:= |s|^{r-1}s$ for
some $r\in (0,1)$, then (\ref{1.1}) reduces back to the classical
fast-diffusion equation (see e.g. \cite{Ar}). For more general
stochastic evolution equations in Hilbert space we refer to
\cite{DRRW,Kim,KR,RRW}  and the references within.

In general, for a fixed number $r\in (0,1)$, we assume that there exist locally
bounded positive measurable functions $\dd,\eta: [0,\infty)\to \R^+$ such that
\beq\label{1.2} \beg{split}
&|\Psi(t,s)|\le \eta_t
(1+|s|^r),\ \ \ s\in\R, t\ge 0,\\
&2\big(\Psi(t,s_1)-\Psi(t,s_2)\big)(s_1-s_2)\ge
\dd_t|s_1-s_2|^2(|s_1|\lor |s_2|)^{r-1},\ \  s_1,s_2\in \R, t\ge
0,
\end{split}\end{equation}
where $\dd$ satisfies $\inf\limits_{t\in[0,T]}\dd_t>0$ for any $T>0$.
 Due to the mean-valued
theorem and the fact that $r<1$, one has $(s_1-s_2)(s_1^r-s_2^r)\ge
r |s_1-s_2|^2(|s_1|\lor |s_2|)^{r-1}$. So, a  simple example for
(\ref{1.2}) is that $\Psi(t,s)=\frac{\dd_t}{2r}s^r$ with $\eta_t=
\frac{\dd_t}{2r}$.

According to \cite[Theorem 3.9]{RRW},  for any $x\in H$ the equation
(\ref{1.1}) has a unique solution $X_t(x)$ with $X_0(x)=x$, which is
a continuous adapted process on $H$ satisfying
$$\E \sup_{t\in [0,T]}\|X_t(x)\|_H^2<\infty,\ \ \ T>0.$$
We intend to establish Harnack inequalities for
$$P_tF:= \E F(X_t),\  F\in \mathfrak{M}_b(H), \ t>0,$$
where $\mathfrak{M}_b(H)$ is the class of bounded measurable
functions on $H$. Since the dissipativity condition (\ref{1.2})  is
essentially weaker than the corresponding one satisfied by the
porous medium situation where $r>1$, the method and results in
\cite{W07} do not apply automatically in the present case.

As in \cite{W07}, we assume that $Q_t(\oo)$ is non-degenerate for
$t> 0$ and $\oo\in\OO$; that is, $Q_t(\oo) x=0$ implies $x=0$. Let
$$\|x\|_{Q_t}:= \beg{cases} \|y\|_2, &\text{if}\ y\in L^2(\m), Q_t
y=x,\\
\infty, &\text{otherwise.}\end{cases}$$

\beg{thm}\label{T1.1}  Assume $(\ref{1.0})$ and $(\ref{1.2})$ with $r\in (0,1).$  If there exist a
constant $\si\ge \ff 4 {r+1}$ and a strictly positive
 function $\xi\in C([0,\infty))$ such that
\beq\label{1.3}
\|x\|_{r+1}^{2}\cdot \|x\|_H^{\si-2}\ge \xi_t \|x\|_{Q_t}^{\si},\ \
\ x\in L^{r+1}(\m),\ t\ge 0
\end{equation} holds on $\OO$, then  for any $T>0$, $P_T$
is strong Feller and for any positive  $F\in \mathfrak{M}_b(H)$,
    $p>1$  and $x,y\in H$,
 \beq\label{1.4} \beg{split}    \ff{(P_TF)^p(y)}{P_TF^p(x)}\le
 \exp\[\frac{p-1}{4}& \(2\int_0^T\theta_t\d t+\lambda_T
T+
 \|x\|_H^2+\|y\|_H^2\)\\
 & +(p-1)\frac{\int_0^T[(\sigma+2)g_t]^2\d t}{4(\sigma\int_0^Tg_t\d
t)^2}\|x-y\|_H^{2}\\ &+
\lambda_T^{\frac{2-\sigma}{2}}\(\frac{\sigma+2}{\sigma}\)^{\sigma+1}
\frac{[2p(p+1)]^{\sigma/2}}{4(p-1)^{\sigma-1}\(\int_0^Tg_t\d
t\)^\sigma}\|x-y\|_H^{\sigma}
\]
\end{split}
\end{equation}
holds for
$$\lambda_T:=\frac{1}{2}\e^{-\int_0^T(2\gamma_t+2q_t+1)\d t}\inf\limits_{t\in[0,T]}\delta_t,\
\theta_t:=q_t+2^{\frac{r+2}{r}}\eta_t^{\frac{r+1}{r}}\delta_t^{-\frac{1}{r}},\
 g_t:=(\delta_t\xi_t)^{\frac{1}{\sigma}}\e^{-\int_0^t\gamma_s\d s}.$$
\end{thm}

\paragraph{Remark.} (1) The right hand side of (\ref{1.4}) comes from our  argument and calculations,
in particular the coupling method modified from \cite{W07} where the case $r\ge 1$ was studied. Comparing with known Gaussian type 
bounds in finite-dimensions, the first two terms in the exponential are natural. The third term of $\|x-y\|_H^\si$ for $\si\ge 4/(r+1)$ seems more technical, which appears when we handle the exponential moment of an additional term in the coupling by using the dissipasitivity of the drift 
$L\Psi$. It is not clear whether this term is exact or not, but since $\si>2$ it does not destroy the short distance (or short time) behaviors of the heat kernel.   On the other hand, due to the weaker  dissipasitivity of the drift, it is reasonable for the semigroup to have worse long time behaviors. 

(2) Harnack inequalities of type (\ref{1.4}) has many applications. For instance, for diffusions on manifolds it has been applied to study the heat kernel estimate (cf. \cite{GW01}), log-Sobolev inequalities and contractivities of the semigroup (cf. \cite{W97, A98, RW1}), and entropy-transportation inequalities (cf. \cite{BGL}). In the symmetric case it was also applied to the study of transition probability kernels for infinite-dimensional diffusions (cf. \cite{AZ, AK}). Here,  due to the weaker dissipasitivity of the drift and the non-symmetry of the semigroup, we are not be able to derive stronger properties like hypercontractivity, which fails even in finite-dimensions e.g. the semigroup generated by
$\DD- \nn |\cdot|^{r+1}$ on $\R^d$ for 
$r<1,$ which is included in our model when $E$ contains finite many elements. Thus, below we only present an application to moment estimates on heat kernels. 
To this end,  we consider the
following time-homogenous case.

\beg{thm}\label{T1.2} Assume $(\ref{1.0}), (\ref{1.2})$ and that the
embedding $L^{r+1}(\m)\subset H$ is compact. Let $\gg\leq 0$ be
constant and   $\Psi, Q$ be deterministic and time-independent.

$(1)$ The Markov semigroup $P_t$ has a unique  invariant probability
measure
 $\mu$   and\\
$\mu(e^{\varepsilon_0\|\cdot\|_H^{r+1}}+\|\cdot\|_{r+1}^{r+1})<\infty$
for some $\varepsilon_0>0$. If $\gamma<0$ then
$\mu(e^{\varepsilon_0\|\cdot\|_H^{2}})<\infty$ for some
$\varepsilon_0>0$.

$(2)$ If $(\ref{1.3})$ holds for some constant $\xi>0$, then $\mu$
has full support on $H$ and for any $x\in H$, $T>0$ and $p>1$, the
transition density $p_T(x,y)$ of $P_T$ w.r.t $\mu$ satisfies
\beg{equation*}\beg{split} \|p_T(x,\cdot)\|_{L^p(\mu)} \le
\bigg\{\int_H\exp\[& -\frac{1}{4(p-1)}\(2\theta T+\lambda_T
T+\|x\|_H^2+\|y\|_H^2\)\\
 & -\frac{\int_0^T[(\sigma+2)g_t]^2\d t}{ 4(p-1)(\sigma\int_0^Tg_t\d
t)^2}\|x-y\|_H^{2}\\
&-\lambda_T^{\frac{2-\sigma}{2}}\(\frac{\sigma+2}{\sigma}\)^{\sigma+1}
\frac{2^{\ff \sigma 2 -2}[p(2p-1)]^{\sigma/2}}{(p-1)\(\int_0^Tg_t\d
t\)^\sigma}\|x-y\|_H^{\sigma}
\] \mu(\d
y) \bigg\}^{-(p-1)/p}.
\end{split}\end{equation*}
where
$$\lambda_T=\frac{\delta}{2}\e^{-(2\gamma+2q+1)T},\ \
\theta=q+2^{\frac{r+2}{r}}\eta^{\frac{r+1}{r}}\delta^{-\frac{1}{r}},\
\  \
 g_t=(\delta\xi)^{\frac{1}{\sigma}}\e^{-\gamma
 t}.$$
\end{thm}

The above two theorems will be proved in the next section by
modifying the argument in \cite{W07}.  To apply Theorems \ref{T1.1}
and \ref{T1.2}, one has to verify (\ref{1.3}) and the compactness of
the embedding $L^{r+1}(\m)\subset H$. Since $r<1$ so that the norm
in $L^{r+1}(\m)$, which is induced by the first drift term in
(\ref{1.1}), is normally incomparable with that in $H$, for
(\ref{1.3}) and the compactness of the embedding we shall need a
Nash (or Sobolev) inequality. Along this line,  explicit sufficient
conditions for the main results to hold, as well as concrete
examples, are presented in Section 3.

\section{Proofs of Theorems \ref{T1.1} and \ref{T1.2}}

\subsection{Proof of Theorem \ref{T1.1}} As explained in \cite{W07}, to
prove the Harnack inequality for $P_t$, it suffices to construct a
coupling $(X_t, Y_t)$, which is a continuous adapted process on
$H\times H$  such that
\\

\noindent(i) $X_t$ solves (\ref{1.1}) with $X_0=x$;\\
(ii) $Y_t$ solves the equation
$$\d Y_t = \big\{L\Psi(t,Y_t)+ \gg_t
Y_t\big\}\d t +Q_t \d \tilde{W}_t,\ Y_0=y$$ for a cylindrical
Brownian motion $\tilde{W}_t$ on $L^2(\mathbf{m})$ under a weighted
probability measure $R\mathbb{P}$, where $\tilde{W}_t$ as well as
the density $R$ will be constructed later on  by a
Girsanov transformation;\\
(iii) $X_T=Y_T,$  a.s.
\\

 As soon as (i)-(iii) are satisfied, then
\beq\label{h}\begin{split}
P_TF(y)&=\mathbb{E}RF(Y_T)=\mathbb{E}RF(X_T)\\
       &\leq(\mathbb{E}R^{p/(p-1)})^{(p-1)/p}(\mathbb{E}F(X_T)^p)^{1/p}\\
       &=(\mathbb{E}R^{p/(p-1)})^{(p-1)/p}(P_TF^p(x))^{1/p}
\end{split}\end{equation}
which implies the desired Harnack inequality provided
$\mathbb{E}R^{p/(p-1)}<\infty$.

 To realize the above idea, for $\varepsilon>0$ and
$\beta\in \mathbf{C}([0,\infty); \mathbf{R^+})$, let $Y_t$ solve the
equation
 \beq\label{2.1} \d
Y_t=\Big\{L\Psi(t,Y_t)+\gg_tY_t+\frac{\beta_t(X_t-Y_t)}{\xy_H
^\varepsilon}\mathbf{1}_{\{t<\tau\}}\Big\}\d t+Q_t\d W_t , \ Y_0=y,
\end{equation}
where $X_t:=X_t(x)$ and $\tau:=\inf\{t\geq0: X_t=Y_t\}$ is the coupling time.

 According to \cite[Theorem 3.9]{RRW}, we can prove that (\ref{2.1}) also has a unique
strong solution $Y_t(y)$ by using the same argument as in \cite[Theorem A.2]{W07}.
Hence, we have $X_t=Y_t$ for $t\geq \tau$ by the pathwise uniqueness of the
solution.

 Let
$$\zeta_t:=
\ff{\bb_tQ_t^{-1}(X_t-Y_t)}{\|X_t-Y_t\|_H^\vv}\mathbf{1}_{\{t<\tau\}}.$$
We have
$$\d Y_t =(L\Psi(t,Y_t) +\gg_tY_t)\d t + Q_t(\d W_t+\zeta_t\d t),
\ \ \ Y_0=y.$$ According to the Girsanov theorem,
$\tilde{W}_t:=W_t+\int_0^t\zeta_s\d s$ is a cylindrical Brownian
motion under $R\mathbb{P}$, where \beq\label{R}
R:=\exp\[-\int_0^T\langle\zeta_t, \d
W_t\rangle-\frac{1}{2}\int_0^T\|\zeta_t\|_2^2\d t\] .
\end{equation}

So, to verify (ii) and (iii), we need to find out $\varepsilon>0$ and $\beta$ such that\\

$\text{(a)}\  \tau\leq T $  a.s.;

$\text{(b)}\ \mathbb{E}\exp\[\lambda\int_0^T\|\zeta_t\|^2_2\d t\]
<\infty,\ \ \lambda>0 $.\\

By (\ref{1.2}) and the It\^{o} formula (see \cite[Theorem
A.2]{RRW}), \beq\label{2.2}\beg{split}   \d\xy_H^2   \leq \big\{&
-\delta_t\mxy\\
&+2\gamma_t\xy_H^2-2\beta_t\xy_H^{2-\varepsilon}\mathbf{1}_{\{t<\tau\}}\big\}\d
t. \end{split}\end{equation} This implies
 \beq\beg{split}\label{2.3}
&\d\big\{\xy_H^2\e^{-2\int_0^t\gg_s\d s}\big\}\\
&\leq -\e^{-2\int_0^t\gg_s\d
s}\big\{\delta_t\mxy+2\beta_t\xy_H^{2-\varepsilon}\mathbf{1}_{\{t<\tau\}}\big\}\d
t.
\end{split}\end{equation}

\beg{lem}
 If $\beta$ satisfies  $\int_0^T\beta_t\e^{-\varepsilon\int_0^t\gamma_s\d s}\d t
 \geq\frac{1}{\varepsilon}\|x-y\|_H^\varepsilon$, then $X_T=Y_T
 $ a.s.
 \begin{proof}
 By (\ref{2.2}) and (\ref{2.3}) we have
$$\frac{2}{\varepsilon}\d\big\{ \xy_H^2\e^{-2\int_0^t\gamma_s\d s}\big\}^{\varepsilon/2}
\leq-\beta_t\e^{-\varepsilon\int_0^t\gamma_s\d s}\d t , \ t\leq\tau\wedge T .$$ If
$T<\tau(\omega)$ for some $\omega$, then
$$\|X_T(\omega)-Y_T(\omega)\|_H^\varepsilon \e^{-\varepsilon\int_0^t\gamma_s\d s}-\|x-y\|_H^\varepsilon
\leq-\varepsilon\int_0^T\beta_t \e^{-\varepsilon\int_0^t\gamma_s\d s }\d
t\leq-\|x-y\|_H^\varepsilon.$$
This implies $X_T(\omega)=Y_T(\omega)$, which is
contradictory to $T<\tau(\omega)$.
\end{proof}
\end{lem}

 From now on, we take $\varepsilon=\frac{\sigma}{\sigma+2}$ and
$$\beta_t=c(\varepsilon\delta_t\xi_t)^{1/\sigma}
\e^{-\frac{2}{\sigma+2}\int_0^t\gamma_s\d s},\ \
c=\frac{\|x-y\|_H^\varepsilon}{\varepsilon\int_0^T
 (\varepsilon\delta_t\xi_t)^{\frac{1}{\sigma}}\e^{-\int_0^t\gamma_s\d s}\d t},$$
 so that (a) holds according to Lemma 2.1. Let $f_t:=\mx$.
By (\ref{2.3}), the H\"{o}lder inequality and (\ref{1.3})  we obtain
\ce \d\big\{\xy_H^2\e^{-2\int_0^t\gamma_s\d s}\big\}^\varepsilon
&\leq&-\varepsilon\delta_t \e^{-2\varepsilon\int_0^t\gamma_s\d s}
\xy_H^{2(\varepsilon-1)}
\mxy \d t\\
&\leq&-\varepsilon\delta_t \e^{-2\varepsilon\int_0^t\gamma_s\d s}
\xy_H^{2(\varepsilon-1)}
\frac{\xy_{r+1}^2}{\mx} \d t\\
 &\leq&-\varepsilon\delta_t\xi_t \e^{-2\varepsilon\int_0^t\gamma_s\d s}
 \frac{\xy_{Q_t}^\sigma}{\xy_H^{\sigma-2\varepsilon}f_t}\d t\\
&=&-\varepsilon\delta_t\xi_t \e^{-2\varepsilon\int_0^t\gamma_s\d s}
\frac{\xy_{Q_t}^\sigma}{\xy_H^{\sigma\varepsilon}f_t}\d t\\
&=&-\frac{\beta_t^\sigma\xy_{Q_t}^\sigma}{c^\sigma\xy_H^{\sigma\varepsilon}f_t}\d
t. \de
Combining this with the H\"{o}lder inequality and the fact that

$$\sup_{a>0}\{a^{\ff{\si-2}\si}b^{\ff 2 \si} -a\}= \Big[\Big(\ff{\si -2}\si\Big)^{\ff{\si-2}2}- \Big(\ff{\si -2}\si\Big)^{\ff \si 2}\Big]b\le b,\ \  b>0$$
implies  $a^{\ff{\si-2}\si}b^{\ff 2 \si}\le a + b, \ a,b>0,$
 we arrive at
 \beq\label{ga}\beg{split}
\int_0^T\|\zeta_t\|_2^2\d t&=\int_0^T\frac{\beta_t^2\xy_{Q_t}^2}{\xy_H^{2\varepsilon}}\d t\\
&\leq\(\int_0^Tf_t^{\frac{2}{\sigma-2}}\d
t\)^{\frac{\sigma-2}{\sigma}}
\(\int_0^T\frac{\beta_t^\sigma\xy_{Q_t}^\sigma}{\xy_H^{\sigma\varepsilon}f_t}\d t\)^{\frac{2}{\sigma}}\\
&\leq\(\int_0^Tf_t^{\frac{2}{\sigma-2}}\d
t\)^{\frac{\sigma-2}{\sigma}}
\(c^\sigma\|x-y\|_H^{2\varepsilon}\)^{\frac{2}{\sigma}}\\
&\leq \lambda\int_0^Tf_t^{\frac{2}{\sigma-2}}\d t+\lambda^{(2-\sigma)/2}c^\sigma\|x-y\|_H^{2\varepsilon},\ \ \lambda>0.
\end{split}
\end{equation}
Since $\sigma\geq\frac{4}{1+r}$ implies $\frac{2}{\sigma-2}\leq\frac{1+r}{1-r}$,  we have

$$f_t^{\ff2 {\si -2}}\le \m\big(1+|X_t|^{r+1}\lor |Y_t|^{r+1}\big)^{\ff{ 2(1-r)}{(\si-2)(1+r)}}
\le \m\big(1+|X_t|^{r+1}\lor |Y_t|^{r+1}\big).$$ Thus,
 \beq\label{2.4}\beg{split}
&\mathbb{E}\exp\[\lambda \int_0^Tf_t^{\frac{2}{\sigma-2}}\d t\]\\
&\leq\mathbb{E}\exp\[\lambda
\int_0^T(1+\|X_t\|_{r+1}^{r+1}+\|Y_t\|_{r+1}^{r+1})\d t\],\ \ \ \lambda>0.
\end{split}\end{equation}
Therefore, to verify (b) we need to prove that
$\int_0^T(\|X_t\|_{r+1}^{r+1}+\|Y_t\|_{r+1}^{r+1})\d t$ is
exponentially integrable. This follows from the following Lemma.

\beg{lem} We have \beq\label{2.5}\mathbb{E}\exp\[\lambda_T
\int_0^T\|X_t\|_{r+1}^{r+1}\d t\]\leq \exp\[\int_0^T\theta_t\d
t+\|x\|_H^2\],
\end{equation}
\beq\beg{split}\label{2.6}
&\mathbb{E}\exp\[\lambda_T \int_0^T\|Y_t\|_{r+1}^{r+1}\d t\]\\
&\leq \exp\[\int_0^T\theta_t\d
t+\|y\|_H^2+\|x-y\|_H^{2(1-\varepsilon)}\int_0^T\beta_t^2
e^{-2\varepsilon\int_0^t\gamma_s\d s}\d t\],
\end{split}\end{equation}
where $\lambda_T=\frac{1}{2}\exp[-\int_0^T(2\gamma_s+2q_s+1)\d
s]\inf\limits_{t\in[0,T]}\delta_t$.

\end{lem}

\beg{proof} Since assumption $(\ref{1.2})$ implies

\beg{equation*}\beg{split} -2\<\Psi(t, X_t),X_t\> & = -2 \<\Psi(t, X_t)- \Psi(t,0), X_t-0\>- 2 \Psi(t,0)\m(X_t)\\
&\le - \dd_t \|X_t\|_{r+1}^{r+1} + 2\eta_t \|X_t\|_{r+1},\end{split}\end{equation*}
by the It\^{o} formula we obtain
\beq\beg{split}\label{ito} \d \|X_t\|_H^2&\leq
\big\{-\delta_t\|X_t\|_{r+1}^{r+1}+2\eta_t\|X_t\|_{r+1}+2\gamma_t\|X_t\|_H^2+q_t\big\}\d
t
+2\<X_t, Q_t\d W_t\>\\
&\leq\Big\{\theta_t-\frac{\delta_t}{2}\|X_t\|_{r+1}^{r+1}+2\gamma_t\|X_t\|_H^2\Big\}\d
t +2\<X_t, Q_t\d W_t\>.
 \end{split}\end{equation}
 Recall that
$\theta_t=q_t+2^{\frac{r+2}{r}}\eta_t^{\frac{r+1}{r}}\delta_t^{-\frac{1}{r}}
$ and $q_t\geq \sup\limits_{\omega\in\Omega}\|Q_t\|_{L_{HS}}^2$.
This implies
$$\begin{aligned} &\d\big\{\e^{-\int_0^t(2\gamma_s+2q_s)\d s}\x_H^2\big\}\\
&\leq \e^{-\int_0^t(2\gamma_s+2q_s)\d s}
\Big\{\theta_t-\frac{\delta_t}{2}\|X_t\|_{r+1}^{r+1}-2q_t\|X_t\|_H^2\Big\}\d
t+2\e^{-\int_0^t(2\gamma_s+2q_s)\d s}\<X_t, Q_t\d W_t\>.
\end{aligned}$$ Hence,
$$\beg{aligned}
&\frac{\delta_{0,T}}{2}\e^{-\int_0^T(2\gamma_s+2q_s)\d s}\int_0^T\x_{r+1}^{r+1}\d t\\
&\leq \int_0^T\theta_t\d
t+\|x\|_H^2+M_T-\int_0^T2q_t\e^{-\int_0^t(2\gamma_s+2q_s)\d
s}\|X_t\|_H^2\d t
\end{aligned}$$
where $\delta_{0,T}:=\inf\limits_{t\in[0,T]}\delta_t$ and $
M_T=2\int_0^T\e^{-\int_0^t(2\gamma_s+2q_s)\d s}\langle X_t, Q_t\d
W_t\rangle$. It is easy to check from (\ref{ito}) and (\ref{1.0})
that $M_t$ is a martingale.
 Taking
$\lambda_T=\frac{\delta_{0,T}}{2}\e^{-\int_0^T(2\gamma_s+2q_s+1)\d
s}$, we obtain
$$\beg{aligned}&\mathbb{E}\exp\[\lambda_T \int_0^T\|X_t\|_{r+1}^{r+1}\d t\]\\
&\leq \exp\[\int_0^T\theta_t\d t+\|x\|_H^2\]
\mathbb{E}\exp\[M_T-\int_0^T2q_t\e^{-\int_0^t(2\gamma_s+2q_s)\d
s}\|X_t\|_H^2\d t\].
\end{aligned}$$
Since $\langle M\rangle_t\leq
\int_0^T4q_t\e^{-\int_0^t(2\gamma_s+2q_s)\d s}\|X_t\|_H^2\d t$ and
$\mathbb{E}\exp[M_t-\frac{1}{2}\langle M\rangle_t]=1$, we obtain
$$ \mathbb{E}\exp\[M_T-\int_0^T2q_t\e^{-\int_0^t(2\gamma_s+2q_s)\d
s}\|X_t\|_H^2\d t\]\leq 1. $$
Thus, (\ref{2.5}) holds.


 Similarly, since (\ref{2.2}) implies $\|X_t-Y_t\|_H^2\leq \e^{2\int_0^t\gamma_s\d
s}\|x-y\|_H^2$, by (\ref{2.1}) and the It\^{o} formula we have
 \ce
&&\d\big\{\e^{-\int_0^t(2\gamma_s+2q_s+1)\d s}\y_H^2\big\}\\
&\leq& \e^{-\int_0^t(2\gamma_s+2q_s+1)\d s}
\[\theta_t-\frac{\delta_t}{2}\|Y_t\|_{r+1}^{r+1}-(2q_t+1)\|Y_t\|_H^2\\
& &+2\y_H\beta_t\xy_H^{1-\varepsilon}\mathbf{1}_{\{t<\tau\}}\]\d t+\d M^\prime_t\\
&\leq&\e^{-\int_0^t(2\gamma_s+2q_s+1)\d s}
\[\theta_t-\frac{\delta_t}{2}\|Y_t\|_{r+1}^{r+1}-2q_t\|Y_t\|_H^2
+\beta_t^2\xy_H^{2(1-\varepsilon)}\mathbf{1}_{\{t<\tau\}}\]\d t+\d M^\prime_t\\
&\leq&\e^{-\int_0^t(2\gamma_s+2q_s+1)\d s}
\[\theta_t-\frac{\delta_t}{2}\|Y_t\|_{r+1}^{r+1}-2q_t\|Y_t\|_H^2
+\beta_t^2 \e^{2(1-\varepsilon)\int_0^t\gamma_s\d s}
\|x-y\|_H^{2(1-\varepsilon)}\]\d t+\d M^\prime_t,
 \de
 where
$M^\prime_t:=\int_0^t2\e^{-\int_0^s(2\gamma_u+2q_u+1)\d u}\<Y_s, B_s
\d W_s\>$ is a martingale. This implies

\beg{equation*}\beg{split}
&\frac{\delta_{0,T}}{2}\e^{-\int_0^T(2\gamma_s+2q_s+1)\d s}\int_0^T\y_{r+1}^{r+1}\d t\\
& \leq \int_0^T\theta_t\d
t+\|y\|_H^2+\|x-y\|_H^{2(1-\varepsilon)}\int_0^T\beta_t^2\e^{-2\varepsilon\int_0^t\gamma_s\d
s}\d t \\
&\qquad+M^\prime_T-\int_0^T2q_t\e^{-\int_0^t(2\gamma_s +2q_s+1)\d
s}\|Y_t\|_H^2\d t.\end{split}\end{equation*} Therefore, taking
$\lambda_T=\frac{\delta_{0,T}}{2}\e^{-\int_0^T(2\gamma_s+2q_s+1)\d
s}$ and noting that

$$\langle M^\prime\rangle_T\leq
\int_0^T4q_t\e^{-\int_0^t(2\gamma_s+2q_s+1)\d s}\|Y_t\|_H^2\d t,$$
we obtain (\ref{2.6}).
\end{proof}

Now, combining (\ref{h}) and (\ref{R}) we obtain
\beq\label{2.8}\beg{split}
& \frac{(P_TF(y))^p}{P_TF^p(x)}\leq \(\mathbb{E}R^{p/(p-1)}\)^{p-1}\\
& =\left\{\mathbb{E}\exp\[\frac{p}{p-1}\int_0^T\langle\zeta_t,
\d W_t\rangle-\frac{p}{2(p-1)}\int_0^T\|\zeta_t\|_2^2\d t\]\right\}^{p-1}\\
& \leq \left\{\mathbb{E}\exp\[ \frac{qp}{p-1}\int_0^T\langle
\zeta_t, \d
W_t\rangle-\frac{q^2p^2}{2(p-1)^2}\int_0^T\|\zeta_t\|_2^2\d t\]
\right\}^{\frac{p-1}{q}}\\
& \ \cdot\left\{\mathbb{E}\exp\[\frac{qp(qp-p+1)}{2(q-1)(p-1)^2}\int_0^T\|\zeta_t\|_2^2\d t\]\right\}^{\frac{(q-1)(p-1)}{q}}\\
& =
\left\{\mathbb{E}\exp\[\frac{qp(qp-p+1)}{2(q-1)(p-1)^2}\int_0^T\|\zeta_t\|_2^2\d
t\]\right\}^{\frac{(q-1)(p-1)}{q}},\ \  q>1.
\end{split}\end{equation}
Moreover, letting
$\lambda=\frac{\lambda_T(q-1)(p-1)^2}{pq(pq-p+1)}$, by (\ref{ga}),
 (\ref{2.4}) and Lemma 2.2 we obtain that \beq\beg{split}
 &
 \mathbb{E}\exp\[\frac{qp(qp-p+1)}{2(q-1)(p-1)^2}\int_0^T\|\zeta_t\|_2^2\d t\]\\
 & \leq \mathbb{E}\exp\[  \frac{\lambda_T}{2}\int_0^T(1+\|X_t\|_{r+1}^{r+1}+\|Y_t\|_{r+1}^{r+1})\d t\\
 & +\frac{qp(qp-p+1)}{2(q-1)(p-1)^2}\(\frac{\lambda_T(q-1)(p-1)^2}{pq(pq-p+1)}\)^{\frac{2-\sigma}{2}}
c^\sigma\|x-y\|_H^{2\varepsilon}\]\\
& \leq \exp\[ \frac{1}{2}\(2\int_0^T\theta_t\d
t+\lambda_TT+\|x\|_{H}^{2}+\|y\|_{H}^{2}+\|x-y\|_H^{2(1-\varepsilon)}\int_0^T\beta_t^2
e^{-2\varepsilon\int_0^t\gamma_s\d s}\d t\)\\
 & \qquad\qquad+\frac{qp(qp-p+1)}{2(q-1)(p-1)^2}\(\frac{\lambda_T(q-1)(p-1)^2}{pq(pq-p+1)}\)^{\frac{2-\sigma}{2}}
c^\sigma\|x-y\|_H^{2\varepsilon}\].
\end{split}\end{equation}
Combing this with (\ref{2.8}) and simply letting $q=2$, we obtain
\beq\beg{split} & \frac{(P_TF(y))^p}{P_TF^p(x)} \leq\exp\[
\frac{p-1}{4}\(2\int_0^T\theta_t\d
t+\lambda_TT+\|x\|_{H}^{2}+\|y\|_{H}^{2}+\\
& \|x-y\|_H^{2(1-\varepsilon)}\int_0^T\beta_t^2
e^{-2\varepsilon\int_0^t\gamma_s\d s}\d t\)
 +\frac{p(p+1)}{2(p-1)}\(\frac{\lambda_T(p-1)^2}{2(p+1)}\)^{\frac{2-\sigma}{2}}
c^\sigma\|x-y\|_H^{2\varepsilon}\]. \end{split}\end{equation} Then
the desired result (\ref{1.4}) follows by using the definition of
$\beta_t$ and $c$.

Finally, since
$$
 |P_TF(y)-P_TF(x)|=|\mathbb{E}(R-1)F(X_T)| \le \|F\|_\infty
\mathbb{E}|R-1|,
$$
and since due to (\ref{2.8}) $R$ is uniformly integrable for
$\|x-y\|_H\leq1$, by the dominated convergence theorem we obtain

$$\lim_{y\to x}|P_TF(y)-P_TF(x)|\leq\|F\|_\infty\lim_{y\to x} \mathbb{E}|R-1|
=\|F\|_\infty\mathbb{E}\lim_{y\to x}|R-1|=0$$ for any bounded
measurable function $F$ on $H$, where the last equality follows from \\
$\lim\limits_{y\to x}R=1$ due to (\ref{ga}). So, $P_T$ is strong
Feller. Now the proof is complete.  \qed

\subsection{Proof of Theorem \ref{T1.2}}

Since the embedding $L^{r+1}(\m)\subset H$ is compact, the existence
and uniqueness of the invariant measure $\mu$ follow
 from b) of \S 4 in \cite{DRRW} for $\m(|\cdot|^{r+1})$ in place of
 $\m(|\cdot|^2).$  So, for (1) it suffices to prove the desired concentration property.
 By (\ref{1.2}) we have
  \ce
\d\x_H^2&\leq&(c-\theta\x_{r+1}^{r+1}+2\gamma\x_H^2)\d t+2\<X_t,
Q\d W_t\>\\
&\leq&(c-\theta\x_{r+1}^{r+1})\d t+2\<X_t, Q\d W_t\>,
 \de
  where $c, \theta>0$ are two constants.
  Then, by a standard argument as in \cite{W07} we obtain
$$\mu(\|\cdot\|_{r+1}^{r+1})<\infty .$$
 If $\gamma<0$ and $\varepsilon_0$ is small enough, then the It\^{o} formula
implies \ce \d \e^{\varepsilon_0\x_H^2}&\leq&\Big(c
-\theta\x_{r+1}^{r+1}+2\gamma\x_H^2+\frac{\varepsilon_0}{2}q\x_H^2\Big)\varepsilon_0
e^{\varepsilon_0\x_H^2}\d t+\d M_t\\
&\leq&(c_1-\theta_1e^{\varepsilon_0\x_H^2})\d t+\d M_t \de for some
constants $c_1, \theta_1>0$ and local martingale
$M_t:=2\varepsilon_0\int_0^t \e^{\varepsilon_0\|X_s\|_H^2}\langle
X_s, Q\d W_s\rangle$.
 This implies
 $$\mu_n(\e^{\varepsilon_0\|\cdot\|_H^2})=\frac{1}{n}\int_0^n\mathbb{E}\e^{\varepsilon_0\|X_t(0)\|_H^2}\d t
 \leq\frac{c_1}{\theta_1}$$
 for $\mu_n:=\frac{1}{n}\int_0^n(\delta_0P_t)\d t$.
  Since $\mu$ is
 the weak limit of a subsequence of $\mu_n$ (see the proof of \cite[Proposition 2.2]{RRW}),
 we obtain  $\mu(\e^{\varepsilon_0\|\cdot\|_H^{2}})<\infty$.  Similarly, $\mu(e^{\varepsilon_0\|\cdot\|_H^{r+1}})<\infty$ holds
 for
 $\gamma=0$.

Finally, the full support of $\mu$ and the $L^p$-estimate of the transition density
follows from the Harnack inequality (\ref{1.4})
 by repeating the proof of Theorem 1.2 (1) and (2) in \cite{W07}.
 \qed

\section{Explicit sufficient conditions and Examples}

To provide explicit sufficient conditions for (\ref{1.3}), we need
the following Nash  inequality:
 \beq\label{N}
\|f\|_2^{2+4/d}\leq C\langle f, -Lf\rangle, \ f\in\mathcal{D}(L),\
\mathbf{m}(|f|)=1.
 \end{equation}

 \beg{lem}\label{L1.2}
 Let $r\in(0, 1)$. Assume that $(\ref{N})$ holds for some
 $d\in(0,\frac{2(r+1)}{1-r})$ and $-(-L)^{1/n}$ is a Dirichlet operator
 for some $n\geq1$. Then the embedding $L^{r+1}(\m)\subset H$ is
 compact. In particular,
$$\|x\|_H:=\langle x, (-L)^{-1}x\rangle^{1/2} \leq c\|x\|_{r+1},
\ x\in L^{r+1}(\mathbf{m})$$
 holds for some $c>0$.  \end{lem}

\begin{proof}  Take  $\vv\in (0,1)$ such that $d_\vv:=
d/\vv \in (d, \frac{2(r+1)}{1-r})$, and let $L_\vv:= - (-L)^\vv.$
By \cite[Theorem 1.3]{BM} and (\ref{N}),

$$\|f\|_2^{2+4/d_\vv}\leq C'\langle f, -L_\vv f\rangle, \ f\in\mathcal{D}(L_\vv),\
\mathbf{m}(|f|)=1$$ holds for some constant $C'>0.$ By this and
again \cite[Theorem 1.3]{BM} and (\ref{N}),
 we have \beq\label{N1} \|f\|_2^{2+\frac{4}{d_\vv n}}\leq c_0\langle f,
(-L_\vv)^{1/n}f\rangle, \ \ f\in\mathcal{D}((-L_\vv)^{1/n}),\
\mathbf{m}(|f|)=1
 \end{equation}
 for some $c_0>0$. Let $T_t$ be the  semigroup generated by $-(-L_\vv)^{1/n}$, which is  sub-Markovian since
$-(-L_\vv)^{1/n}= - (-L)^{\vv/n}$ is a Dirichlet operator. Then it
follows from (\ref{N1}) that (see \cite{Davies})
$$\|T_t\|_{1\rightarrow\infty}\leq c_1t^{-d_\vv n/2}, \ \ t>0$$
holds for some constant $c_1>0$. Since $T_t$ is contractive in $L^1(\m)$ and  symmetric  in $L^2(\m)$, by this and the Riesz-Thorin interpolation theorem we obtain 

\beq\label{OO} \|T_t\|_{1\to 2} = \|T_t\|_{2\to\infty} \le c_2 t^{-d_\vv n/4},\ \ \ t>0\end{equation} for some constant $c_2>0.$ Moreover,
since $\lambda_1>0$ so that $\|T_t\|_{2\to 2} \le \e^{-\ll_1 t}$ for $t>0,$ (\ref{OO}) yields 
$$\|T_t\|_{1\rightarrow\infty}\leq \|T_{t/4}\|_{1\rightarrow 2}\|T_{t/2}\|_{2\rightarrow 2}
\|T_{t/4}\|_{2\rightarrow\infty}\leq c_3t^{-d_\vv
n/2}e^{-\lambda_1^{\frac{\vv}{n}}t/2}, \ \ t>0$$ for some $c_3>0.$ By this and the
Riesz-Thorin interpolation theorem we conclude that for any $1<p<q$,
$$\|T_t\|_{p\rightarrow q}\leq\|T_t\|_{1\rightarrow\infty}^{\frac{q-p}{pq}}
\leq c_4[t^{-d_\vv
n/2}e^{-\lambda_1^{\frac{\vv}{n}}t/2}]^{\frac{q-p}{pq}}, \ \ t>0$$
holds for some constant $c_4>0$. Therefore,
$$C_{p,q}:=\int_0^\infty\|T_t\|_{p\rightarrow q}dt<\infty $$
provided $\frac{q-p}{pq}<\frac{2}{d_\vv n}$. Thus,
$$\|(-L_\vv)^{-1/n}\|_{p\rightarrow q}\leq C_{p,q}<\infty,\ \ \ \frac{q-p}{pq}<\frac{2}{d_\vv n}.$$
Since $d_\vv<\frac{2(r+1)}{1-r}$, letting
$p_i:=\frac{r+1}{1-2(i-1)(r+1)/d_\vv n}\ (1\leq i\leq n+1)$, one has
$$p_1=r+1,\ \ \frac{p_{i+1}-p_i}{p_{i+1}p_i}=\frac{2}{d_\vv n}\ (1\leq i\leq n)\ \ \text{and}
\ \  p_{n+1}=\frac{r+1}{1-2(r+1)/d_\vv}>\frac{r+1}{r}. $$ So, there
exist
$r+1=:p_1^\prime<p_2^\prime<\cdots<p_{n+1}^\prime:=\frac{r+1}{r}$
such that $\frac{p_{i+1}^\prime-p_i^\prime}{p_{i+1}^\prime
p_i^\prime}<\frac{2}{d_\vv n},\ 1\leq i\leq n.$ Therefore,
$$c^2:=\|(-L_\vv)^{-1}\|_{r+1\rightarrow (r+1)/r}\leq
\prod_{i=1}^n \|(-L_\vv)^{-\frac{1}{n}}\|_{p_i^\prime\rightarrow
p_{i+1}^\prime}\leq \prod_{i=1}^n
C_{p_i^\prime,p_{i+1}^\prime}<\infty.$$ This implies
$$\begin{aligned} &
\langle
x,(-L_\vv)^{-1}x\rangle\leq\|x\|_{r+1}\|(-L_\vv)^{-1}x\|_{(r+1)/r}\\
&\leq\|x\|_{r+1}^2\|(-L_\vv)^{-1}\|_{r+1\rightarrow
(r+1)/r}=c^2\|x\|_{r+1}^2,\ \ x\in L^{r+1}(\mathbf{m}).
\end{aligned}$$ Then the proof is completed since $\{x\in L^2(\m):\  \<x,
(-L_\vv)^{-1}x\>\le N\}$ is relatively compact in $H$ for any $N>0$.
\end{proof}

\beg{cor}\label{C1.2}
 Let $Q_te_i=q_ie_i$ for $i\geq1$ with
$\sum_{i=1}^\infty\frac{q_i^2}{\lambda_i}<\infty$, so that Q is
Hilbert-Schmidt from $L^2(\mathbf{m})$ to H. If
$\varepsilon\in(0,1)$ and L satisfies $(\ref{N})$ for some $d\in(0
,\frac{2\varepsilon(1+r)}{1-r})$, $-(-L)^{1/n}$ is a Dirichlet
operator
 for some $n\geq1$ and there exist $c>0$,
$\sigma\geq\frac{4}{1+r}$ such that
$$q_i\geq c\lambda_i^{\frac{\sigma+2\varepsilon-2}{2\sigma}}, \ i\geq1 .$$
Then the embedding $L^{r+1}(\m)\subset H$ is compact and
$(\ref{1.3})$ holds for the same $\sigma$.
\end{cor}

\begin{proof} By Lemma \ref{L1.2}, it suffices to verify
(\ref{1.3}). By the H\"{o}lder inequality,
\beq\label{3.1}\beg{split}
\|x\|_Q^\sigma&=\(\sum_{i=1}^\infty\langle
x,e_i\rangle^2q_i^{-2}\)^{\sigma/2}=\(\sum_{i=1}^\infty\frac{\langle
x,e_i\rangle^2}{\lambda_i^{\frac{\sigma-2}{\sigma}}}\lambda_i^{\frac{\sigma-2}{\sigma}}q_i^{-2}\)^{\sigma/2}\\
&\leq\(\sum_{i=1}^\infty\langle
x,e_i\rangle^2\lambda_i^{\frac{\sigma-2}{2}}q_i^{-\sigma}\)\(\sum_{i=1}^\infty\frac{\langle
x,e_i\rangle^2}{\lambda_i}\)^{\frac{\sigma-2}{2}}\\
&=\|x\|_H^{\sigma-2}\(\sum_{i=1}^\infty\langle
x,e_i\rangle^2\lambda_i^{\frac{\sigma-2}{2}}q_i^{-\sigma}\)\\
&\leq c^{-\sigma}\|x\|_H^{\sigma-2}\(\sum_{i=1}^\infty\langle
x,e_i\rangle^2\lambda_i^{-\varepsilon}\) .
 \end{split}\end{equation}
By (\ref{N}) and \cite[Theorem 1.3]{BM},  there exists a constant
$C_\varepsilon>0$ such that
$$\|f\|_2^{2+4\varepsilon/d}\leq C_\varepsilon\langle f, (-L)^\varepsilon f\rangle ,
\ f\in\mathcal{D}((-L)^\varepsilon),\ \mathbf{m}(|f|)=1.$$ Applying
Lemma 3.1 to $-(-L)^\varepsilon$ in place of $L$, there exists a
constant $c_1>0$ such that
$$\|x\|_{r+1}^2\geq c_1\|(-L)^{-\varepsilon/2}x\|_2^2=c_1\sum_{i=1}^\infty\langle x, e_i\rangle^2\lambda_i^{-\varepsilon}
.$$ Combining this with (\ref{3.1}),  we obtain (\ref{1.3}) for some
constant $\xi>0$.

\end{proof}

\beg{exa} Let $Q_te_i=q_ie_i (i\geq1)$, $\Psi(t,x):=|x|^{r-1}x$ and $\gamma_t:=c$
for some constant $c<0$. Let $L:=\Delta$ be the Laplace operator on a bounded domain
in $\mathbf{R}$ with Dirichlet boundary conditions. If $r\in(\frac{1}{3},1)$, then
all assertions in Theorem $1.1$ and $1.2$ hold provided there exist constants $c_1,
c_2>0$, $\alpha<1$ and $\varepsilon\in(\frac{1-r}{2(1+r)}, \frac{r}{1+r})$ such that
$$ c_1i^{\varepsilon(r+1)+1-r}\leq q_i^2\leq c_2i^\alpha ,\ i\geq1 .$$
\end{exa}
\beg{proof}Since $\lambda_i\geq ci^2$ for some constant $c>0$, by
Corollary 3.2 the above conditions imply $(\ref{1.3})$.
\end{proof}

\noindent$\mathbf{Remark}:$ We can also consider the case where $L$
is the Laplace operator on a post critical finite self-similar
fractal with $s>0$ the Hausdorff dimension of the fractal in the
effective resistance metric. In this case we has $\lambda_i\geq
ci^{(s+1)/s},\ i\geq1$ for some $c>0$ according to \cite[Theorem
2.11]{HK}.
\\

To construct examples for our results on high dimensional spaces, we may e.g. take
$L=-(-\Delta)^\alpha$ for large enough $\alpha>0$. More generally, let $-L_0$ be a
self-adjoint Dirichlet operator on $L^2(\m)$ with discrete spectrum
$$(0<) \ll^{(0)}_1\le \ll_2^{(0)}\le \cdots$$
and the corresponding unit eigenfunctions $\{e_i\}_{i\ge 1}$ be an ONB on $L^2(\m)$.
As in Corollary \ref{C1.2}, let $Q e_i:= q_i e_i$ for a sequence $\{q_i\ne 0\}_{i\ge
1}.$ Let, for simplicity, $\gamma_t=-c_0$ and $\Psi\in C(\R)$ satisfy
\beq\label{4.1} \beg{split} &|\Psi(s)|\le \eta
(1+|s|^r),\ \ \ s\in\R,\ t\ge 0,\\
&2\<\Psi(s_1)-\Psi(s_2), s_1-s_2\>\ge \dd\m(|s_1-s_2|^2(|s_1|\lor |s_2|)^{r-1}),\ \
s_1,s_2\in \R,\ t\ge 0
\end{split}\end{equation}
for some
 $c_0\ge 0$ and $\eta,\dd>0$.  For any positive constant $\alpha$,
  we consider the equation (\ref{1.1}) for
$$L:= - (-L_0)^\alpha=-\sum_{i=1}^\infty (\ll_i^{(0)})^\aa
\<e_i,\cdot\> e_i.$$ That is, consider \beq\label{4.2} \d X_t=
-\big\{(-L_0)^\alpha\Psi(X_t) +c_0 X_t\big\}\d t + Q\d W_t.\end{equation}

\beg{prp} \label{P4.1} Let $L_0$ satisfy $(\ref{N})$ with
$d\in(0,\frac{2\varepsilon(1+r)}{1-r})$ for some
$\varepsilon\in(0,1)$, and $\Psi$ satisfy $(\ref{4.1})$. If there
exists $\aa>\frac{d(1-r)}{2\varepsilon(1+r)}$ such that
\beq\label{4.3} \sum_{i=1}^\infty
\ff{q_i^2}{(\ll_i^{(0)})^\aa}<\infty,\ \
  q_i\geq
c(\lambda_i^{(0)})^{\frac{\aa(\sigma+2\varepsilon-2)}{2\sigma}}, \
i\geq1, \end{equation}
where $\sigma\geq\frac{4}{1+r}$ is a
constant, then the Markov semigroup of the solution to $(\ref{4.2})$
satisfies all assertions in Theorems $\ref{T1.1}$ and $\ref{T1.2}$
for the same $\sigma$ and some $\xi>0$.
\end{prp}

\beg{proof} We only need to notice that the eigenvalues of
$L:=-(-L_0)^\aa$ are
 $$-\ll_i:=-(\ll_i^{(0)})^\aa,\ \ i\ge 1.$$
 Obviously, all conditions in Corrollary \ref{C1.2} are satisfied
 for the present situation, hence (\ref{1.3}) holds.
 \end{proof}

To conclude this paper, we present an example where $L_0$ is the
Dirichlet Laplacian on a finite volume domain in $\R^d$, so that $L$
can be taken as high order differential operators on a domain.

\paragraph{Example 3.5.} In the situation of Proposition \ref{P4.1} but simply take $q_i=i^\theta, i\ge 1$
 where $\theta>0$ is a constant. Let $L_0:= \DD$ be the Dirichlet Laplace operator on a domain
$D\subset \R^d$ with finite volume,  and let $\m$ be the normalized
volume measure on $D$. By the Sobolev inequality we have (see
\cite[Corollaries 1.1 and 3.1]{W00})
$$\ll_i^{(0)} \ge c i^{2/d},\ \ \ i\ge 1$$
for some $c>0$. If $\sigma\geq\frac{4}{1+r}$ and
$$\theta >\max\left\{\frac{\sigma+2\varepsilon-2}{4(1-\varepsilon)},
\frac{(\sigma+2\varepsilon-2)(1-r)}{2\sigma\varepsilon(1+r)}\right\}$$
for some $\varepsilon\in(0,1)$,  then (\ref{4.3}) holds for any
$\aa\in\(\frac{(2\theta+1)d}{2}\vee
\frac{(1-r)d}{2\varepsilon(1+r)},\  \frac{\sigma\theta
d}{\sigma+2\varepsilon-2}\]$, so all assertions in Theorems
\ref{T1.1} and \ref{T1.2} hold for the solution to (\ref{4.2})
according to Proposition \ref{P4.1}.

\paragraph{Acknowledgement.} The authors would like to thank the referee for useful comments.

\end{document}